\documentclass{tran-l}
 \newtheorem{thm}{Theorem}[section]
 \newtheorem{cor}[thm]{Corollary}
 \newtheorem{lem}[thm]{Lemma}
 
 \theoremstyle{definition}
 \newtheorem{defn}[thm]{Definition}
 \theoremstyle{remark}
 \newtheorem{rem}[thm]{Remark}
 \numberwithin{equation}{section}
\newcommand{\J}{\mathcal{J}}
\newsymbol\lesssim 132E
\newsymbol\gtrsim 1326
\newsymbol\varnothing 203F
\newsymbol\varsubsetneq 2320
\newsymbol\varsupsetneq 2321
\newcommand{\dR}{I \! \! R}
\def\a{{\textsf{a}}}
\def\esssup{{\rm{ess}}\!\!\!\sup}
\def\dist{{\rm{dist}}}
\def\diam{{\rm{diam}}}
\def\loc{{\!\rm{loc}}}
\def\Q{{\mathcal Q}}
\def\oo{{\mathcal O}}
\def\I{{\mathcal I}}
\def\pij{{\partial_{ij}}}
\def\phk{{\partial_{hk}}}
\def\Lz{{{\mathcal{L}}_0}}
\def\L{{\mathcal{L}}}
\def\Dp{{\mathcal{D}}_p}
\def\Ll{{{\mathcal{L}}_1}}
\def\W{{\wp}}
\def\ww{{\W(x)\over \W(B(x))}}
\def\K{{\mathcal{K}}}
\def\C{{\mathcal{C}}}
\def\O{{\Omega}}
\def\Ol{{O(\lambda, \Lambda)}}
\def\Ole{{O(\lambda, \Lambda, \eta)}}

\def\oz{{\omega_0}}
\def\ol{{\omega_1}}

\def\b{{\partial B}}

\def\Gaq{{{\Gamma_\alpha(Q)}}}
\def\dd{{\partial D}}
\def\bb{{\partial B}}
\def\tz{{\triangle_0}}
\def\tj{{\triangle_j}}
\def\tjz{{\triangle_j^0}}
\def\drq{{\triangle_r(Q)}}
\def\dtq{{\triangle_{2r}(Q)}}
\def\trq{{T_r(Q)}}
\def\ttr{{T_{2r}(Q)}}
\def\vmo{{\rm{V\!M\!O}}}
\def\bmo{{\rm{B\!M\!O}}}
\def\bmor{{\rm{B\!M\!O\!_{\varrho}}}}
\def\bme{{\rm{B\!M\!O\!_{\varrho_0}}}}
\def\bms{{{\rm B\!M\!O\!}_{\varrho^*}}}
\def\bmss{{{\rm B\!M\!O}_{2\varrho^*}}}
\def\r{{\dR^n}}
\begin{document}

\title[The $L^p$ Dirichlet Problem]
 {The $L^p$ Dirichlet Problem and Nondivergence Harmonic Measure}

\author{ Cristian Rios }

\address{School of Mathematics, University of Minnesota, Minneapolis, MN 55455}
\curraddr{Department of Mathematics and Statistics,
McMaster University, Hamilton, ON, L8R-B19 Canada}

\email{riosc@math.mcmaster.ca}



\date{April 3, 2002}
\dedicatory{In memory of E.Fabes}

\copyrightinfo{2002}{American Mathematical Society}%

\begin{abstract}
We consider the Dirichlet problem
$$\left\{
\begin{array}{rcl}
 \L u & = & 0\quad\text{ in }D\cr
u & = & g\quad\text{ on }\dd
\end{array}
\right.$$ for two second order elliptic operators $\L_k u=\sum_{i,j=1}^na_k^{i,j}(x)\,\pij u(x)$, $k=0,1$, in
a bounded Lipschitz domain $D\subset\r$. The coefficients $a_k^{i,j}$ belong to the space of bounded mean oscillation
$\bmo$ with a suitable small $\bmo$ modulus. We assume that $\Lz$ is regular in $L^p(\dd,d\sigma)$ for some $p$,
$1<p<\infty$, that is, $\|Nu\|_{L^p}\le C\,\|g\|_{L^p}$ for all continuous boundary data $g$. Here $\sigma$ is the
surface measure on $\dd$ and $Nu$ is the nontangential maximal operator. The aim of this paper is to establish
sufficient conditions on the difference of the coefficients  $\varepsilon^{i,j}(x)=a^{i,j}_1(x)-a^{i,j}_0(x)$ that will
assure the perturbed operator $\Ll$ to be regular in $L^q(\dd,d\sigma)$ for some $q$, $1<q<\infty$.
\end{abstract}

\maketitle

\section{Introduction}

In the present note we consider linear elliptic second order differential operators in nondivergence form
$\L=\sum_{i,j=1}^na^{i,j}(x)\,\pij$, where $A(x)=(a^{i,j}(x))_{i,j=1}^n$ is a symmetric matrix verifying
the uniform ellipticity and boundedness condition
\begin{equation}\label{ellip}
\lambda\, |\xi|^2\le \xi^t A(x) \xi\le \Lambda\, |\xi|^2,\qquad x,\, \xi\in\r
\end{equation}
for some fixed $0<\lambda\le\Lambda<\infty$ and $n\ge 2$. We study the Dirichlet problem
\begin{equation}\label{Dirich}
\left\{
\begin{array}{rcl}
 \L u & = & 0\quad\text{ in }D\cr
u & = & g\quad\text{ on }\dd
\end{array}
\right.
\end{equation}
on a bounded Lipschitz domain $D\subset\r$. From~\cite{CFL2} and a standard approximation argument it follows that if
the coefficients $a^{i,j}$ are in $\vmo$ ($\bme$) and $g\in\C(\dd)$ problem (\ref{Dirich}) has a unique solution
$u=u_g\in \C(\overline{D})\bigcap W^{2,p}_\loc(D)$ for all $p$, $1<p<\infty$ ($1<p<p_0(\varrho_0)$). 
We denote by $\sigma$ be the surface measure on $\dd$ and we say that the operator $\L$ is regular in
$L^p(\dd,d\sigma)$ or that $\Dp$ holds for $\L$ in $D$, $1<p<\infty$, if there exists a constant $C_p$ which depends on
$n$, $\lambda$, $\Lambda$, $D$, $p$ and the $\bmo$ modulus of the coefficients such that for all continuous boundary
data $g$ the solution $u$ of (\ref{Dirich}) verifies
\begin{equation}\label{DePe}
\| Nu\|_{L^p(\dd,d\sigma)}\le C_p\, \| g\|_{L^p(\dd,d\sigma)},
\end{equation}
where $Nu$ is the nontangential maximal operator
$$Nu(Q)=\sup_{\Gaq}|u(x)|$$
here and henceforth $\Gaq$ denotes the interior truncated cone (of opening $\alpha$)
\begin{equation}\label{gaq}
\Gaq=\{ x\in D:|x-Q|\le(1+\alpha)\,\delta(x)\}\bigcap B_{r^*}(Q),
\end{equation}
$\delta(x)=\dist(x,\dd),$ $B_r(x)$ denotes the ball in $\r$ centered at $x$ of radius $r$ and
$\alpha\ge\alpha^*=\alpha^*(D)>0$, $r^*=r^*(D,\lambda,\Lambda,\eta)>0$ are fixed (here $\eta$ is the $\bmo$ modulus of the
coefficients of $\L$, see Section~\ref{prelim} and (\ref{tcone})). When necessary, we will write $N_\alpha u$ for the
nontangential maximal operator of opening $\alpha$.

The purpose of this note is to give sufficient conditions for the preservation of the regularity of the $L^p$ Dirichlet
problem under small perturbations on the coefficients. Given two elliptic operators
$\L_k=\sum_{i,j=1}^na_k^{i,j}(x)\partial_{x_ix_j}$, where $A_k(x)=(a_k^{i,j}(x))_{i,j=1}^n$, $k=0,1$, are symme\-tric
matrices verifying (\ref{ellip}), let $\varepsilon(x)=(a^{i,j}_1(x)-a^{i,j}_0(x))_{i,j=1}^n$ be the difference between
the coefficients and $B(x)=B_{\delta(x)/2}(x)$, $x\in D$, we consider the quantity
\begin{equation}\label{adef}
\displaystyle{\a(x)=\max_{1\le i,j\le n}\esssup_{y\in B(x)}|\varepsilon^{i,j}(y)}|.
\end{equation}
For $Q\in\dd$ and $r>0$ we denote the boundary ball of radius $r$ at $Q$ by $\drq=B_r(Q)\bigcap\dd$, and the Carleson
region at $Q$ of radius $r$ by $\trq=B_r(Q)\bigcap D$. Our main result is the following:
\begin{thm}\label{main}
Suppose that $\Lz$ verifies $\Dp$ for some $p$, $1<p<\infty$, then there exists
$\varrho_0=\varrho_0(n,\lambda,\Lambda,D,C_p)>0$ such that if $a_k^{i,j}\in\bme(\r)$, $1\le i,j\le n$, $k=0,1$, and
\begin{equation}\label{carl}
\displaystyle{\sup_{Q\in\dd,\, r>0}{1\over\sigma(\trq)}\int_\trq {\a^2(x)\over\delta(x)}\, dx=M<\infty,}
\end{equation}
then $\Ll$ verifies ${\mathcal{D}}_{q}$ for some $q$, $1<q<\infty$.
\end{thm}

A similar result was established in~\cite{FeKP} for divergence form operators with coefficients in $L^\infty(\r)$. We
are able to adapt the divergence case techniques and obtain the results in~\cite{D},~\cite{RF} and~\cite{FeKP} (under
extra assumptions on the coefficients) for the nondivergence case (see~\cite{Cr}).  This gives a partial answer to the problem posed by
C. Kenig in~\cite{K2} (Problem 3.3.9). Condition (\ref{carl}) says that the measure $\a^2/\delta\, dx$ is a
Carleson measure with respect to $\sigma$ with Carleson norm bounded by $M$.

By the maximum principle the correspondence $g\mapsto u_g(x)$ is a positive linear functional on $\C(\dd)$ for each fixed $x\in D$ .
The Riesz representation theorem implies that there exist a unique regular positive Borel measure
$\omega^x=\omega^x_{\L,D}$ such that
$$u(x)=\int_\dd g(Q)\, d\omega^x(Q).$$
The measure $\omega^x$ is called the {\em harmonic measure} for $\L$ and $D$ at
$x$ and constitutes one of our main tools in the proof of Theorem~\ref{main}. Also crucial for this task is the concept of {\em normalized adjoint solution} (n.a.s.),
first introduced in~\cite{Ba} (see also~\cite{Ba2},~\cite{FSt}). In~\cite{EK} n.a.s. are used to define a proper area
function for solutions of nondivergence form operators with bounded coefficients. We also use the theory of
Muckenhoupt weights~\cite{Mu},~\cite{GCRF} and in particular the result in~\cite{Eu} which establishes that non\-ne\-ga\-ti\-ve
adjoint solutions are $A_p$ weights for all $p$, $p_0\le p\le\infty$, where $p_0$ depends on the $\bmo$ modulus of the coefficients. Other important elements in our proofs are the a
priori estimates for solutions~\cite{CFL},~\cite{CFL2}, basic properties of the harmonic measure~\cite{Ba2},~\cite{FGMS},~\cite{FSt} and weighted Poincar\'e inequalities~\cite{FKS}.

\begin{rem} It is known~\cite{D1}~\cite{D2} that the Laplacian operator $\Delta=\sum_{i=1}^n\partial^2_{x_i^2}$
is regular in $L^p$ for $2-\varepsilon<p<\infty$ where $\varepsilon=\varepsilon(n,D)$. On the other hand, examples
in~\cite{MMo} show the existence of a nondivergence operator $\Ll$ with continuous coefficients $A_1$ in the closure of
the unit ball $B$ in $\r$, such that $\Ll=\Delta$ on $\partial B$ and the harmonic measure $\ol=\omega_{\Ll,B}$ is
singular with respect to the surface measure $\sigma$. In particular, $\Ll$ is not regular in $L^p(\dd,d\sigma)$ for
any $p$ (see Theorem~\ref{DP}). Setting $\Lz=\Delta$, the modulus $\a(x)$ corresponding to this example violates condition (\ref{carl}). This shows that the perturbation problem addressed in Theorem~\ref{main} is non trivial, even for continuous coefficients.
\end{rem}


\section{Preliminaries}\label{prelim}

In general, we write $X\lesssim Y$ when there exists a constant $C>0$ which depends at most on $n$,  $\lambda$,
$\Lambda$, $\eta$ and $D$
such that $X\le C\, Y$. Similarly, we define the expression
$X\gtrsim Y$ and write $X\approx Y$ when $X\lesssim Y$ and $X\gtrsim Y$.

If $G\subset\r$ is a Borel set we denote by $\C(G)$ the space of real valued continuous functions on $G$. If $\mu$ is a
$\sigma$-finite Borel measure on $G$, $L^p(G,d\mu)$, $1\le p<\infty$ denotes the Banach space of $\mu$-measurable
functions $f$ on $G$ such that $\| f\|_{L^p(G,d\mu)}=(\int_G|f|^p\, d\mu)^{1\over p}<\infty$. We use $dx$ to denote the
Lebesgue measure in $\r$, $|E|=\int_E dx$ for any Borel set $E$ in $\r$ and we write $L^p(G)=L^p(G,dx)$. The spaces
$L^\infty(G,d\mu)$, $L^p_\loc(G,d\mu)$ are also defined in a standard way. If $G\subset\r$ is open, $k$ is a
nonnegative integer and $1\le p\le\infty$ we set $W^{k,p}(G)$ to be the Sobolev space of functions $f$ with $k$ weak
derivatives in $\L^p(G)$ (see~\cite{GT} Chapter 7).

Given $f\in\L^1_\loc(\r)$ we set
$$\eta(r,x)=\eta_f(r,x)=\sup_{s\le r}{1\over |B_s(x)|}\int_{B_s(x)}|f(y)-f_{B_s(x)}|\,dy$$
where $f_{E}={1\over |E|}\int_{E}f(y)\,dy$. We say that $f$ has {\em bounded mean oscillation} or that $f\in\bmo(\r)$
if $\eta\in L^\infty(\dR^+,\r)$ and set $\|f\|_{\bmo(\r)}=\|\eta_f\|_{L^\infty((\dR^+,\r))}$.

\begin{defn}\label{etaset} Given $\varrho>0$ and $\zeta>0$, we let $\Phi(\varrho,\zeta)$ be the set
$$\Phi(\varrho,\zeta)=\{\eta:\dR^+\mapsto \dR^+,\, \eta\text{ non-decreasing, }\eta(r)\le\varrho\text{ whenever
}r<\zeta\}.$$ We also set $\Phi(\varrho)=\bigcup_{\zeta>0}\Phi(\varrho,\zeta)$, and given $\eta\in\Phi(\varrho)$ we
denote by $\zeta(\eta,\varrho)=\zeta(\eta)=\sup\{\zeta>0:\eta\in\Phi(\varrho,\zeta)\}$.
\end{defn}

If $\varrho >0$ we say that $f\in\bmor(\r)$ if $\eta(r)=\|\eta(r,\cdot)\|_{L^\infty(\r)}$ lies
in $\Phi(\varrho)$. If $\lim_{r\rightarrow 0^+}\eta(r)=0$ we say that $f$ has {\em vanishing
mean oscillation} or that $f\in\vmo(\r)$ (see~\cite{Sr}). We also define $\bmo(G)$ and $\bmo(G,d\mu)$ in a standard way through the modulus
$$\eta(r,x,G,\mu)=\sup_{s\le r}{1\over\mu(B_s(x)\bigcap G)}\int_{B_s(x)\bigcap G}|f(y)-f_{B_s(x)\bigcap G,d\mu}|\,d\mu,$$
where $f_{E,d\mu}={1\over \mu(E)}\int_{E}f(y)\,d\mu$, $G\subset\r$ is a Borel set and $\mu$ is a Borel measure .

Given an non-decreasing function $\eta: \overline{\dR^+}\mapsto\overline{\dR^+}$, we denote by $\Ole$ the class of
operators ${\L}=\sum_{i,j=1}^n a^{i,j}(x)\,\partial_{x_ix_j}$, w    ith sy\-mme\-tric coefficients $A(x)=(a^{i,j}(x))_{i,j=1}^n$
verifying the ellipticity and boundedness conditions (\ref{ellip}) and such that $a^{i,j}\in\bmo(\r)$, $1\le i,j\le n$,
with $\bmo$-modulus of continuity $\eta$ in $D$. When there is no restriction on the regularity of the coefficients of
$\L$, we say $\L\in\Ol$.

We denote by $D$ a bounded Lipschitz domain in $\r$. That is, a bounded, connected open set $D$ such that its boundary
$\dd$ can be covered by a finite number of open right circular cylinders whose bases have positive distance from $\dd$
and corresponding to each cylinder $C$ there is a coordinate system $(x',x_n)$ with $x'\in\dR^{n-1}$, $x_n\in\dR$ with
$x_n$ axis parallel to the axis of $C$, and a function $\psi:\dR^{n-1}\mapsto \dR$ satisfying a Lipschitz condition
($|\psi(x')-\psi(y')|\le m_0\, |x'-y'|$) such that $C\bigcap D=\{ (x',x_n): x_n>\psi(x')\}\bigcap C$, and
$C\bigcap\dd=\{ (x',x_n):x_n=\psi(x')\}\bigcap C$. Whenever we say that a quantity depends on $D$, we mean it depends
on the Lipschitz character of $D$. In what follows we assume that $D$ is contained in the unit ball and contains the
origin.

Let $\triangle$ denote a generic boundary ball in $\dd$, i.e. $\triangle=\drq$ for some $r>0$, $Q\in\dd$. Given two
Borel measures $\mu$ and $\nu$ on $\dd$, we say that $\mu$ is in $A_\infty$ with respect to $\nu$ on $\dd$ and we write
$\mu\in A_\infty(d\nu)$ if there exist $0<\zeta<1$ and $\kappa>0$ such that
\begin{equation}\label{ainf}
{\nu(E)\over\nu(\triangle)}>\zeta\Rightarrow {\mu(E)\over\mu(\triangle)}>\kappa,
\end{equation}
whenever $E\subset\triangle$ and $E$ is a Borel set. The theory of $A_\infty$ weights originates in~\cite{Mu} and~\cite{CoFe} where the results below can be found (see also~\cite{GCRF} and~\cite{St2}). We say that $\mu$ is in the
reverse H\"older class $B_{p'}(d\nu)$, $1<p'<\infty$, if $\mu$ is absolutely continuous with respect to $\nu$ and
$k={d\mu\over d\nu}$ verifies
$$\left\{{1\over\nu(\triangle)}\int_\triangle k^{p'}\, d\nu\right\}^{1\over p'}\le C\,
{1\over\nu(\triangle)}\int_\triangle k\, d\nu,$$ for all boundary balls $\triangle\subset\dd$. The weight $k$ is in
$A_p(d\nu)$, $1<p<\infty$ if
\begin{equation}\label{apcond}
\displaystyle{\left\{{1\over\nu(\triangle)}\int_\triangle k\, d\nu\right\} \left\{{1\over\nu(\triangle)}\int_\triangle
k^{-{1\over p-1}}\, d\nu\right\}^{1\over p-1}\le C<\infty.}\end{equation} It is easy to see that $A_\infty$ is an
equivalence relation, and that $k\in A_p(d\nu)$ if and only if $k^{-1}\in B_{p'}(d\mu)$, ${1\over p}+{1\over p'}=1$.
The best constant $C$ in (\ref{apcond}) is called the $A_p(d\nu)$ ``norm'' of $k$ and we denoted it by $
|[k]|_{A_p(d\nu)}$ or $|[\mu]|_{A_p(d\nu)}$. We will also use the convention $k\in A_p$ (resp.: $B_{p'}$, $A_\infty$)
whenever $k\in A_p(d\sigma)$ (resp.: $B_{p'}(d\sigma)$, $A_\infty(d\sigma)$).

We say that a measure $\nu$ is  a {\em doubling} measure, with doubling cons\-tant $c=c(\nu)$ if
$\nu(\triangle_{2r}(Q))\le c\, \nu(\drq)$ for all $r>0$ and $Q\in\dd$. It is also well known that if $\mu\in A_p(d\nu)$
then $\mu$ is a doubling measure if and only if $\nu$ is a doubling measure and $c(\mu)=c(\nu)^p|[\mu]|_{A_p(d\nu)}$.

Given a Borel measure $\mu$ on $\dd$, we denote by $M_\mu g(Q)$ the Hardy-Littlewood maximal operator at $Q$ with respect to $\mu$, that is:
\begin{equation}\label{maxw}M_\mu g(Q)=\sup_{Q\in\triangle(Q)}{1\over\mu(\triangle(Q))}\int_\triangle g(P)\, d\mu(P)
\end{equation}
where $\triangle(Q)$ denotes a generic boundary ball in $\dd$ centered at $Q$.  It is known that if $\mu$ is a doubling measure, then
$$\| M_\mu f\|_{L^p(\dd,d\nu)}\le C_p\, \| f\|_{L^p(\dd,d\nu)},\qquad 1<p\le\infty,$$
with $C_p>0$ independent of $f$, if and only if $\nu\in A_p(d\mu)$~\cite{CoFe} (see also~\cite{GCRF}).

If $u$ is the solution of (\ref{Dirich}) with boundary data $g\in \C(\dd)$ then $Nu\approx M_\omega g$ (\cite{Ba} Theorem 7.3). Here and henceforth $\omega$ denotes the harmonic measure for $\L$ and $D$ at a fixed point $x_0\in D$. 
Since the harmonic measure is a doubling measure~\cite{Ba} (see also~\cite{FSt}) we have that the maximal operator is
bounded in $L^p(\dd,d\omega)$, $1<p\le\infty$, and then $\|Nu\|_{L^p(\dd,d\omega)}\lesssim\|g\|_{L^p(\dd,d\omega)}$ for
all $p$, $1<p\le\infty$. From the weighted maximal theorem~(\cite{GCRF}IV.2.1) we then have that $\L$ verifies $\Dp$
for some $1<p<\infty$ if and only if $\omega$ is a weight in the reverse H\"older class $B_{p'}(d\sigma)$, ${1\over
p}+{1\over p'}=1$.

Other basic fact of the theory of weights is that
$$A_\infty(d\nu)=\bigcup_{p>1} A_p(d\nu)=\bigcup_{p'>1} B_{q'}(d\nu),$$
hence, to prove Theorem~\ref{main} it is enough to show $\ol\in A_\infty(d\sigma)$. The following theorem is a
consequence of the weighted maximal theorem, the theory of weights, and the inequalities $Nu\approx M_\omega g$.

\begin{thm}\label{DP}
Let $\omega$ be the harmonic measure with respect to $\L$ in $D$ and $\mu$ be a Borel measure on $\dd$. The following
are equivalent:
\begin{enumerate}
  \item[{\rm (i) }] $\omega\in A_\infty(d\mu).$
  \item[{\rm (ii) }] There exist $1<p<\infty$ such that $\Dp(d\mu)$ holds, that is
  \[\| Nu\|_{L^p(\dd,d\mu)}\le C_p\, \| g\|_{L^p(\dd,d\mu)}.\]
  \item[{\rm (iii) }] $\omega$ is absolutely continuous with respect to $\mu$ and $k={d\omega\over
d\mu}$ belongs to $B_{q}(d\mu)$, $({1\over p}+{1\over q}=1)$.
\end{enumerate}
\end{thm}


\subsection{A Priori Estimates and Properties of Solutions}
\begin{thm}[Maximum principle~\cite{GT} 9.1.]\label{MP} Let $\L\in\Ol$, $D$ be a bounded domain and
$u\in C(\overline{D})\bigcap W^{2,n}_\loc(D)$ verifies $\L u\ge f$ with $f\in L^n(D)$, then there exists $C>0$ which depends
only on $n$, $\diam(D)$, $\lambda$ and $\Lambda$ such that
$$\sup_D u\le \sup_\dd u^+ + C\, \| f\|_{L^n(D)}.$$
\end{thm}



\begin{thm}
{\label{strong}} Let $w\in A_p$, $p\in(1,\infty)$. There exist
positive numbers $c=c(n,p,\lambda,\Lambda,|[w]|_{A_p})$ and $\tilde{\varrho}_p=\tilde{\varrho}_p(n,c)$, such that if
$\eta\in\Phi(\tilde{\varrho}_p)$, and $\L\in\Ole$, then for any open set $\Omega\subset\r$,
$\diam(\Omega)\le\zeta(\eta)$, and any $u\in W^{2,p}_0(\Omega)$ we have
$$\| \pij u\|_{L^p(\Omega,w)}\le c\, \| \L u\|_{L^p(\Omega,w)}\quad \forall i,j=1,\cdots, n.$$
\end{thm}

\begin{proof} This theorem is an immediate consequence of the techniques in~\cite{CFL} and weighted estimates for
singular integral operators and commutators. We give a sketch of the proof.
Given $u\in W^{2,p}_0(\O)$, we have the following representation formula~\cite{CFL}
\begin{equation}\label{repr}
\begin{array}{rcl}
\displaystyle{\pij u(x)}& = &\displaystyle{\K_{i,j}\left( \sum_{h,k=1}^n (a^{hk}(x)-a^{hk}(\cdot))\,
\phk u(\cdot)+\L u(\cdot)\right)} \cr & & \displaystyle{+ \L u (x)\int_{|t|=1}\Gamma_i(x,t)t_j d\sigma(t).}
\end{array}
\end{equation}
where
$$\displaystyle{\K_{i,j}f(x)} = \displaystyle{\lim_{\varepsilon\rightarrow 0}\int_{\O\backslash B_\varepsilon(x)}
\Gamma_{i,j}(x,x-y)f(y)\, dy}$$ is a principal value operator and $\Gamma_i(x,t)={\partial\over\partial
t_i}\Gamma(x,t)$, $\Gamma_{i,j}(x,t)={\partial^2\over\partial t_i\partial t_j}\Gamma(x,t)$. Here $\Gamma(x,t)$, $x\in
\O$, is a fundamental solution of $\L_0 u(t) = \sum_{i,j=1}^n a^{i,j}(x)\,\pij u(t)$ (see~\cite{CFL} for
details).
For each pair $(i,j)$, $1\le i,j\le n$, $\K_{i,j}$ is a singular integral operator with a regular kernel and then
$\K_{i,j}$ is bounded in $L^p(dw)$ with operator norm which depends on $n$, $\lambda$, $\Lambda$, $p$ and $|[w]|_{A_p}$
(c.f.~\cite{GCRF}, Theorem IV3.1., see also~\cite{CFL} Theorem 2.11).
From the weighted estimates for commutators in~\cite{BrCe}, we have that the commutators $\C_{i,j,h,k}$, $1\le
i,j,h,k\le n$, given by $\C_{i,j,h,k}f(x)=a^{h,k}(x)\K_{i,j}f(x)-\K_{i,j}(a^{i,j}f)(x)$ are bounded in $L^p(\r,dw)$
with norm $\| C_{i,j,h,k}\|\le \tilde{c}\, \| a^{i,j}\|_{\bmo(\r,dw)}$. Moreover, whenever $f$ is supported in $\O$ we
have the localized estimate (see~\cite{CFL} Theorem 2.13)
$$\|\C_{i,j,h,k}f\|_{L^p(\O,dw)}\le c\, \| a^{i,j}\|_{\bmo(\O)}\,\|f\|_{L^p(\O,dw)},$$
where we used that since $w\in A_p$ we have $\tilde{c}\, \| a^{i,j}\|_{\bmo(\O,dw)}\le c\, \|
a^{i,j}\|_{\bmo(\O)}$.
Finally, it is not difficult to check that the factor multiplying $\L u(x)$ in (\ref{repr}) is uniformly bounded, with
bound depending only on $n$, $\lambda$ and $\Lambda$. From (\ref{repr}) and the mentioned estimates we have
$$\| \pij u\|_{L^p(\O,dw)}\le c\, \| \L
u\|_{L^p(\O,dw)}+c\, \sum_{h,k=1}^n \| a^{h,k}\|_{\bmo(\O)}\,\| \phk u \|_{L^p(\Omega,dw)},$$ the theorem follows
taking $\tilde{\varrho}_p\le (2\, n^2\, c)^{-1}$.
\end{proof}

The following theorem follows from the results in~\cite{CFL2}, the techniques just exposed and standard arguments (see Theorem 8.1 in~\cite{Ba}).
\begin{thm}
{\label{strrong}} Let $w\in A_p$, $p\in[n,\infty)$ and $D\subset\r$ be a Lipschitz domain. There exist a positive $\varrho_p=\varrho_p(n,p,\lambda,\Lambda,|[w]|_{A_p})$, such that if $\eta\in\Phi(\varrho_p)$, and $\L\in\Ole$, then for any $f\in L^p(D,w)$, there exists a unique $u\in C(\overline{D})\bigcap W^{2,p}_\loc(D,w)$ such that $\L u=f$ in $D$ and $u=0$ on $\dd$. 

Moreover, if $\dd$ is of class $C^2$, then  $u\in W^{1,p}_{\!0}(D,w)\bigcap W^{2,p}(D,w)$ and there exists a positive $c=c(n,p,\lambda,\Lambda,|[w]|_{A_p},\chi,\dd)$, with $\chi=\diam(D)/\zeta(\eta,\varrho_p)$, such that
$$\| u\|_{W^{2,p}(D,w)}\le c\, \| f\|_{L^p(D,w)}.$$
\end{thm}

For each $x\in D$ and $f$, $u$ as in Theorem~\ref{strrong}, the maximum principle (Theorem~\ref{MP})
implies that the positive linear functional $f\longmapsto -u(x)$ is bounded on $L^p(D)$. From Riesz' representation
theorem we have that there exist a unique nonnegative function $G_{\L,D}(x,\cdot)\in L^{p'}(D)$ with $p'={p\over{p-1}}$
such that
\begin{equation}
u(x)=-\int_D G_{\L,D}(x,y)\, f(y)\, dy.
\end{equation}

\begin{defn}[Green's function]
\label{green} {\rm The function $G_{\L,D}(x,y)$ is called the {\em Green's function} for $\L$ in $D$.  For simplicity we will often write $G(x,y)=G_{\L,D}(x,y)$.}
\end{defn}

\begin{cor}\label{grep} For all $\varphi\in C_{\! c}^\infty(D)$ and $x\in D$ we have
$$\varphi(x)=-\int_D \L\varphi(y)\, G(x,y)\, dy.$$
\end{cor}

\subsection{Properties of the Harmonic Measure}
\begin{lem}[~\cite{Ba},~\cite{FGMS}]\label{lots}Let $\triangle=\drq$, $Q\in\dd$:
\begin{enumerate}
  \item $\triangle'=\triangle_s(Q_0)\subset\triangle$, $x\in D\backslash\ttr$. Then
 $\omega^{x_r(Q)}(\triangle')\approx{\omega^x(\triangle')\over\omega^x(\triangle)}.$
  \item $\omega^{x_r(Q)}(\triangle)\approx 1$, $x_r(Q)$ verifies $\delta(x_r(Q))\approx |x_r(Q))-Q|$.
  \item {\em (Doubling property)}
\label{doubl} $\omega^x(\triangle)\approx\omega^x(\dtq)$, $x\in D\backslash\ttr$.  In particular, the harmonic measure $\omega$ can not have atoms.
\end{enumerate}
\end{lem}


\subsection{Adjoint Solutions and Area Functions}
\begin{defn}[Adjoint solution]\label{adjdef} Given $\L\in\Ol$, a locally integrable function $v$ is an adjoint
solution of $\L$ in a domain $D$ and we write $\L^*v=0$ in $D$, if
$$\int_D v\,\L\varphi\, dx=0$$
for all $\varphi\in C_{\! c}^\infty (D)$. More generally, if $f\in L^1_\loc(D)$, we say that $v$ is a solution to to
$\L^* v=f$ if
$$\int_D v\,\L\varphi\, dx=\int_D f\,\varphi\, dx$$
for all $\varphi\in C_{\! c}^\infty (D)$.
\end{defn}
So $v$ is an adjoint solution for $\L$ in $D$ if $\pij(a^{ij}v)=0$ in the sense of distributions. Suppose now that
$\eta\in\Phi(\varrho_n)$, where $\varrho_n$ is given by Theorem~\ref{strrong} for $p=n$ and $\Phi(\varrho_n)$ is as in
Definition~\ref{etaset}. If $\L\in\Ole$ and $G$ is the Green's function for $\L$ in $D$ (see Definition~\ref{green}),
then for each $x\in D$, $G(x,\cdot)$ is an adjoint solution of $\L$ in $D\backslash\{ x\}$. Indeed, from
Corollary~\ref{grep} we have
$$\int_D\L\varphi(y)\, G(x,y)\, dy=-\varphi(x)=0\qquad\forall\varphi\in C_{\! c}^\infty(D\backslash\{ x\}).$$

The following existence theorem for adjoint solutions is a consequence of the classical theory for smooth
operators~\cite{GT}, Theorem~\ref{strrong} and the maximum principle, we omit the standard proof.

\begin{thm}\label{exadj} Let
$\eta\in\Phi(\varrho_n)$, where $\varrho_n$ is given by Theorem~\ref{strrong} for $p=n$ and $\Phi(\varrho_n)$ is as in
Definition~\ref{etaset}. If $\L\in\Ole$ and $G$ is the Green's function for $\L$ in $D$ then for any $f\in L^{n\over
n-1}(D)$, there exists a solution $v\in L^{n\over n-1}(D)$ to the problem
\begin{equation}\label{adjequ}
\displaystyle{\L^* v =f\quad\text{ in }D.}
\end{equation}
\end{thm}

Examples in~\cite{Ba3} show that even if the coefficients of $\L$ are continuous, adjoint solutions could be not in $L^\infty(D)$. Further ``weight type'' regularity exists in the case of positive nonnegative solutions. In~\cite{FSt} it was
shown that if $w$ is a nonnegative adjoint solution for $\L\in \Ole$, then $\log w$ lies in $\bmo$. Next theorem~\cite{Eu} is a
more precise version of this result, more suitable to our applications.
\begin{thm}{\label{apbq}}
Let $0<\varrho\le\varrho_n$ where $\varrho_n$ is as in Theorem~\ref{strrong}, $\eta\in\Phi(\varrho)$ (see
Definition~\ref{etaset}), ${\L}\in\Ole$, and $w$ be a nonnegative adjoint solution to ${\L}^*w=\pij (a^{i,j}w)=0$ on $B_{10}$. Then there exists $\varrho_0=\varrho_0(n,\lambda,\Lambda,\varrho)>0$, such that
$\log w$ is a function lying in $\bme$. Moreover, $\varrho_0\lesssim\varrho^{\gamma}$ for some
$\gamma=\gamma(n,\lambda,\Lambda)>0$.
\end{thm}
Recall that (the Lipschitz domain) $D\subset B_1$, where $B_r=B_r(0)$.
We pick a point $\overline{x}\in\partial B_9$ and we let $\W=\W(\L)$ be given by
\begin{equation}\label{W}
\W(y) =G_{\L,B_{10}}(\bar{x},y)\quad\text{ in }B_{10}\end{equation}
where $G_{\L,B_{10}}$ is the Green's function for $\L$ in $B_{10}$ (see Definition~\ref{green}).  From the previous theorem we have that there exists $\varrho^{**}>0$ such that if $\L\in\Ole$ with
$\eta\in\Phi(\varrho^{**})$, then $\W$ is a weight in $A_{4\over
3}(B_{8})$, and $|[\W]|_{A_{4\over 3}}$ (see (\ref{apcond})) depend only on $n,\, \lambda,\,\Lambda$.
We set
\begin{equation}\label{rhost}
\varrho^*={1\over 2}\,\min\{ \varrho^{**},\,\varrho_n\},
\end{equation}
where $\varrho_n$ is given by Theorem~\ref{strrong} for $w=\W$ and $p=n$. Since
$|[\W]|_{A_{n}}\le|[\W]|_{A_{2}}\le|[\W]|_{A_{4\over 3}}$, the constant $\varrho^*$ depends only on $n,\, \lambda$ and
$\Lambda$. Note also that $\W\,dx$ is a doubling measure in $B_8$ with doubling constant which depends only on $n,\,
\lambda$ and $\Lambda$.

\begin{defn}[\cite{Ba}, n.a.s.]\label{nas} Let $\L\in\Ol$, a {\em normalized adjoint solution} for $\L^*$
in $D$ is any function $w$ of the form
$$w(x)={v(x)\over \W(x)}$$
where $v$ is a solution of the adjoint equation $\L^*v=0$ in $D$ and $\W$ is given by (\ref{W}).
\end{defn}

Normalized adjoint solutions, first introduced in~\cite{Ba}, enjoy many desirable properties adjoint solutions fail to
verify. Fo\-llo\-wing the techniques in~\cite{Ba}, the Dirichlet problem for n.a.s. is uniquely solvable for continuous
boundary data and coefficients in $\Ole$, with $\eta\in\Phi(\varrho_n)$. A Harnack principle holds for nonnegative
n.a.s., as well as a boundary Harnack inequality and a comparison principle (see~\cite{Ba},~\cite{Ba2}
and~\cite{FGMS}). Although the definition of n.a.s. depends on the particular choice of the {\em normalizing function}
$\W$, this choice has no qualitative impact in our applications.

\begin{lem}\label{cappo} Let $\varrho^*$ be given by (\ref{rhost}), then if $\eta\in\Phi(\varrho^*)$, $\L\in\Ole$, $u$
sa\-tis\-fy $\L u=0$ in $D$, and $v\in L^1_{\loc}(D)$ is a nonnegative adjoint solution for $\L$ in $B_{3r}(x_0)\subset
D$, then for $0<2\,r\le\zeta(\eta)$ and any constants $\beta$ and $\gamma$ the following holds
$$\begin{array}{rcl}
\displaystyle{\int_{B_r(x_0)} |\nabla^2u(x)|^2\, v(x) dx}&\lesssim&r^{-4}\,\displaystyle{\int_{B_{2r}(x_0)}
|u(x)-\beta-\gamma\, x |^2\, v(x) dx}\cr\cr &&\displaystyle{+\,r^{-2}\,\int_{B_{2r}(x_0)} |\nabla(u(x)-\gamma\, x)|^2 \,
v(x) dx.}
\end{array}$$
\end{lem}

\begin{proof} We choose a nonnegative $\phi\in C_{\! c}^\infty(B_{2r}(x_0))$ such that $\phi\equiv 1$ in $B_r(x_0)$ and
$|\partial_j\phi|\le M\, r^{-j}$, $j=0,1,2$, with $M>0$ a universal constant. Applying Theorem~\ref{strong} to the function $(u-\beta-\gamma\,
x)\,\phi$, we have
\begin{equation}\label{ekinter}
\int_{B_r(x_0)} |\nabla^2u(x)|^2 \W(x)\, dx\lesssim\int |\L ((u-\beta-\gamma\, x)\,\phi(x))|^2 \W(x)\, dx.
\end{equation}
Developing the derivatives, re-arranging terms, applying H\"older inequality and since $u$ is a solution for $\L$, we
have that the right hand side of (\ref{ekinter}) is bounded by
$$C\,\int \left( |u(x)-\beta-\gamma\, x|^2 |\nabla^2\phi(x)|^2 + |\nabla(u(x)-\gamma\, x)\cdot \nabla\phi(x)|^2\right) \, \W(x) dx,$$
which proves the lemma in the case $v=\W$. Let now $w=v/\W$, with $v$ as in the statement of Lemma~\ref{cappo}. Then
$w$ is a {\em normalized adjoint solution} of the operator $\L$ in $D$. The lemma follows from Harnack inequality for
n.a.s. (c.f.~\cite{Ba},~\cite{FGMS}). \end{proof}

\begin{lem}[\cite{EK}, Lemma 2]\label{lem2} Let $G(x,y)$ be the Green's function in $D$ for $\L\in\Ol$. Then there
is a constant $r_0$ depending on the Lipschitz character of $D$, such that for all $Q\in\dd$, $r\le r_0$, $y\in\partial
B_r(Q)\bigcap\Gamma_1(Q)$, and $x\not\in T_{4r}(Q)$, the following holds
$${G(x,y)\over\delta(y)^2}{\W(B(y))\over \W(y)}\sim\omega^x(\Delta_r(Q)).$$
\end{lem}
The following lemma establishes that the regularity of the Dirichlet problem depends on the coefficients locally at the boundary.
\begin{lem}\label{wlog} Let $\Lz$, $\Ll\in\Ole$ be such that if $A_0$ and $A_1$ denote their respective matrices of
coefficients, we have that for some $s_0>0$
$$A_1(x)=A_2(x)\quad\text{ for all }x\in D\text{ such that }\delta(x)\le s_0,$$
then there exists $C>0$, depending only on $n$, ellipticity, $D$ and $s_0$ such that
$$C^{-1}\,\oz(\drq)\le \ol(\drq)\le C\,\oz(\drq),\quad\forall Q\in\dd,\, r>0,$$
where $\oz$ and $\ol$ denote the harmonic measures for $\Lz$ and $\Ll$, respectively. In particular,
$\ol\in\bigcap_{p>1} A_p(d\oz)\cap\bigcap_{q>1}B_{q}(d\oz)$.
\end{lem}

\begin{proof} Let $Q\in\dd$, $r={s_0\over 8}$, $\O_s=\{x\in D\, :\, \delta(x)<s\}$, and $G_0$, $G_1$ denote the Green's
functions in $D$ for $\Lz$ and $\Ll$, respectively. Then for any $x\in D$ we have that the functions $G_0(x,\cdot)$ and
$G_1(x,\cdot)$ are adjoint solutions for $\Lz$ in $\O_{s_0}\backslash\{ x\}$.  From the comparison principle for
normalized adjoint solutions~\cite{FGMS}, we have that
$${G_0(x,y)\over G_1(x,y)}\approx{G_0(x,y_r(Q))\over G_1(x,y_r(Q))},\qquad y\in\trq,\ x\in D\backslash\O_{4r}.$$
From Lemma~\ref{lem2}~\cite{EK} we have
$${\oz^x(\triangle_y)\over\ol^x(\triangle_y)}\approx
{\oz^x(\drq)\over\ol^x(\drq)}.$$ where
$\triangle_y=\triangle_s(P)$ with $s\approx\delta(y)\approx |y-P|$, $P\in\dd$.  From Lemma~\ref{lots}-(2) and the
interior Harnack inequality we have 
$$\oz^x(\triangle_y)\approx\ol^x(\triangle_y),\qquad y\in\trq ,\ x\in T_{5r}(Q)\backslash\O_{4r}.$$
From the interior Harnack inequality for solutions, ~\cite{KS}, the harmonic measures $\omega_i^{x}$, $i=0,\, 1$, where
$x\in\O_{5r}\backslash\O_{4r}$, are comparable to the respective harmonic measures at the center of $D$, $\oz$ and
$\ol$ (with constants depending on $D$ and $s_0$); thus we obtain for some $C>0$ as wanted
$$C^{-1}\,\oz(\triangle_y)\le \ol(\triangle_y)\le C\,\oz(\triangle_y),\qquad\forall y\in\O_{s_0\over 8}.$$
The general case follows from Lemma~\ref{lots}-(3). \end{proof}

\begin{defn}[Area functions] For a function $u$ defined on $D$, the area function of aperture $\alpha$, $S_\alpha
u$ and the second area function of aperture $\alpha$, $A_\alpha u$, are defined respectively as

$$S_\alpha u (Q)^2=\int_\Gaq {\delta(x)^2\over \W(B(x))}|\nabla u(x)|^2\, \W(x)\, dx,$$
and
$$A_\alpha u (Q)^2=\int_\Gaq {\delta(x)^4\over \W(B(x))}|\nabla^2u(x)|^2\, \W(x)\, dx$$
where $\W$ is as in {\rm (\ref{W})}, $B(x)= B_{\delta(x)/2}(x)$,  and $Q\in\dd$.
\end{defn}

\begin{thm}[\cite{EK}]\label{ntaf}
Let $\L\in\Ol$, $u\in C(\overline{D})$ a solution to $\L u = 0$ in $D$, $u(Q)=g(Q)$ on $\dd$, where $g\in C(\dd)$. If
$\nu$ is a positive  Borel measure on $\dd$, which is in $A_\infty(\omega)$, where $\omega$ is the harmonic measure for
$\L$ in $D$ evaluated at $0$, then given $0<p<\infty$, $\alpha>0$, $\beta>0$, there exists a constant $C$ which
depends on $n$, ellipticity, $p$, $\alpha$, $\beta$, the $A_\infty$ constant of $\nu$ and the Lipschitz character of
$D$ such that
$$\| S_\alpha u\|_{L^p(\dd, d\nu)} \le C\, \| N_\beta u \|_{L^p(\dd, d\nu)}.$$
Moreover, if $u(0)= 0$
$$\| N_\alpha u\|_{L^p(\dd, d\nu)} \le C\, \| S_\beta u \|_{L^p(\dd, d\nu)}.$$
\end{thm}

\begin{lem}\label{capp2}
Let $\L\in\Ole$, $u$ a solution to $\L u = 0$ in $D$ and $w$ an $A_2$ weight. Then if $B_{2r}(x)\subset D$ and
$2r\le\zeta(\eta)$, we have
$$ \int_{B_r(x)} |\nabla^2 u(y)|^2\, w(y)\, dy \le {C\over r^2}\int_{B_{2r}(x)} |\nabla u(y)|^2\, w(y)\, dy $$
where $C$ depends only on $n$, D, ellipticity and $|[w]|_{A_2}$.
\end{lem}

\begin{proof} Choose a cut-off function $\phi\in C_{\! c}^\infty(B_{2r}(x))$, with $\phi=1$ on $B_r(x)$,
$$\| \partial_i\phi\|_{L^\infty(B_{2r}(x))}\le M\, r^{-i},\qquad i=0,1,2,$$
and $M$ a universal constant. Applying Lemma~\ref{cappo} to $u\,\phi$ we obtain
$$\int_{B_r(x)}|\nabla^2u|^2\, dw\lesssim {1\over r^2}\,\int_{B_{2r}(x)}|\nabla u|^2\, dw +{1\over
r^4}\,\int_{B_{2r}(x)}|u-u_{{B_{2r}(x)},w}|^2 \, dw, $$ where $u_{{B},w}={1\over w(B)}\int_B u\, dw$. By the weighted
Poincar\'e type inequality (Theorem 1.5 in~\cite{FKS}) we have

$$\int_{B_{2r}(x)}|u-u_{{B_{2r}(x)},w}|^2 \, w(y)\, dy \le C\, r^2\, \int_{B_{2r}(x)}|\nabla u|^2 \, w(y)\, dy$$
where $C$ has the required dependence, this finishes the proof. \end{proof}

We will fix from now on the length of the truncation $r^*$ for the cones $\Gaq=\{ x\in
D:|x-Q|\le(1+\alpha)\,\delta(x)\}\bigcap B_{r^*}(Q)$ defined in (\ref{gaq}). We set
\begin{equation}\label{tcone} r^*=\min\{r_0,\,(2\sqrt{n})^{-1}\,\zeta(\eta)\},
\end{equation}
where $r_0=r_0(D,\lambda,\Lambda)$ is as in Lemma~\ref{lem2}, $\eta$ is the common modulus of continuity for $\Lz$ and
$\Ll$ in Theorem~\ref{main} and $\zeta(\eta)$ is given by Definition~\ref{etaset}.
\begin{thm}\label{whit}
Let $\varrho^*$ be given by (\ref{rhost}), $\eta\in\Phi(\varrho^*)$,  $\L\in\Ole$ and $u$ be a solution to $\L u = 0$
in $D$. For any $\alpha\ge\alpha^*(D)>0$ we have the point-wise inequality on $\dd$
$$A_\alpha u(Q)\le C\, S_{c\alpha} u(Q)$$
where $c>1$ depends only on  dimension $n$ and $C$ depends on $n$ and the ellipticity constants.
\end{thm}

\begin{proof} Let $\{Q_j\} _{j=1}^\infty$ be a Whitney decomposition of $D$ into cubes, that is, $D=\bigcup_{j=1}^\infty Q_j$, $Q_j$ and $Q_k$ have disjoint interiors for $j\neq k$, and $r_j=\text{side-length}(Q_j)\approx\dist(Q_j,\dd)$. Denote by $x_j$ the center of $Q_j$. We assume that $B_{2\sqrt{n}r_j}(x_j)\subset D$, and denote $\tilde{Q_j}$ the cube with center $x_j$ and side length $\sqrt{n} r_j$. Note that there exist constants $N$ and $c>1$ depending only on the
dimension $n$ such that
\begin{equation}\label{chin}
\sum_{\stackrel{\displaystyle{j=1}}{Q_j\bigcap\Gaq\neq\emptyset}}^\infty\chi_{\tilde{Q}_j}\le N\,
\chi_{\Gamma_{c\alpha}(Q)}\end{equation} where $\chi_E(x)$ denotes the characteristic function of the set $E$. We have
$$\begin{array}{rcl}
A_\alpha^2 u(Q)&=&\displaystyle{\sum_{Q_j\bigcap\Gaq\neq\emptyset}\int_{Q_j\bigcap\Gaq}{\delta(x)^4\over
\W(B(x))}|\nabla^2u(x)|^2\, \W(x)\, dx}\cr &&\cr
 &\le & \displaystyle{C\, \sum_{Q_j\bigcap\Gaq\neq\emptyset}{r_j^4\over \W(B(x_j))}\int_{Q_j} |\nabla^2 u(x)|^2\, \W(x). dx}
\end{array}$$
Since $\eta\in\Phi(\varrho^*)$, $\W$ is an $A_2$ weight with $A_2$ constant depending only on $n$, $\lambda$ and
$\Lambda$. Note that by assumption (\ref{tcone}) we have $2\sqrt{n}\, r_j\le\zeta(\eta)$, we apply Lemma~\ref{capp2}
and (\ref{chin}) to obtain

$$\begin{array}{rcl}
A_\alpha^2 u(Q)& \le&\displaystyle{ C\, \sum_{Q_j\bigcap\Gaq\neq\emptyset}{r_j^2\over \W(B(x_j))}\int_{\tilde{Q}_j} 
|\nabla u(x)|^2\, \W(x)\, dx} \cr &&\cr & \le&\displaystyle{C\,
\sum_{{Q}_j\bigcap\Gaq\neq\emptyset}\int_{\tilde{Q}_j}{\delta(x)^2\over \W(B(x))}|\nabla u(x)|^2\, \W(x)\, dx}\cr &&\cr &
\le&\displaystyle{N\, C\, \int_{\Gamma_{c\alpha}(Q)}{\delta(x)^2\over \W(B(x))}|\nabla u(x)|^2\, \W(x)\, dx}\cr &&\cr & = &
\displaystyle{N\, C\, S^2_{c\alpha}u(Q).}
\end{array}$$
\end{proof}

We will also find useful an averaged version of the nontangential maximal func\-tion.
\begin{defn}\label{enie}
For a function $u$ in $D$ and $\alpha>0$, define the {\em modified nontangential maximal function} $N_\alpha^0(u)$ by
$$N^0_\alpha u(Q)=\sup_\Gaq\left\{\int_{B_0(x)}u^2(y)\, {\W(y)\over \W(B(y))}dy\right\}^{1\over 2},$$
where $B_0(x)=B_{\delta(x)\over 6}(x)$.
\end{defn}

From the doubling property of the weight $\W$, we have $N^0_\alpha u\lesssim\, N_\alpha u$. When $u$ is a solution of
{\em any} elliptic operator in $\Ole$, we also have:
\begin{lem}\label{enenie} Let $\varrho^*$ be given by {\rm (\ref{rhost})}, $\eta\in\Phi(\varrho^*)$,  $\Ll\in\Ole$,
and $u$ be a solution of $\Ll u=0$ in $D$, then
$$N_\alpha u(Q)\lesssim N_\alpha^0 u(Q)\qquad \forall Q\in\dd.$$
\end{lem}

\begin{proof} Let $Q\in\dd$ and $x\in\Gaq$, from a known reverse H\"older inequality for solutions (\cite{S}, see also~\cite{GT} Theorem 8.17) we have $|u(x)|\lesssim \| u\|_{L^{3/2}(B_0(x))}\,|B_0(x)|^{-3/2}$. Then from H\"older inequality we get
$$\begin{array}{rcl}
|u(x)|
&\lesssim&\displaystyle{{1\over |B_0(x)|^{2\over 3}}\,\left( \int_{B_0(x)}|u|^2\,{\W(y)\over\W(B(y))}\, dy\right)^{1\over 2}\,\left(\int_{B_0(x)}{\W(B(y))^3\over\W(y)^3}\, dy\right)^{1\over 6}  . }
\end{array}$$
Since $\W$ is an $A_{4\over 3}$-weight, we have $\| 1/\W\|_{L^3(B_0(x))}\lesssim |B_0(x)|^{4\over 3}\, \W(B_0(x))^{-1}$,
and from the doubling property of $\W$ we have $\W(B(y))\lesssim\W(B(x))$ for all $y\in B_0(x)$, thus
$$|u(x)|\lesssim \displaystyle{ N_\alpha^0 u(Q)\,{\W(B(x))^{1\over 2}\over |B_0(x)|^{1\over 2}}\, {|B_0(x)|^{1\over 2}\over\W(B_0(x))^{1\over 2}}}\lesssim N_\alpha^0 u(Q).$$
The lemma follows taking supremum for $x\in\Gaq$ on the above inequality. \end{proof}

\section{The main local estimate}\label{proof271}
We present here a version of Theorem~\ref{main} in which surface measure $d\sigma$ is replaced by the harmonic measure
$d\oz$ and (\ref{carl}) is modified accordingly. A similar result was proved in~\cite{FeKP} for divergence form
operators.
\begin{thm}\label{271}Let $G_0(x,y)$ be the Green's function for $\Lz$ in $D$ and set $G_0(y)=G_0(0,y)$.
There exist $\varepsilon_0>0$ which depend only on $n$, $\lambda$, $\Lambda$ and $D$ such that if $a_k^{i,j}\in\bms$,
$1\le i,j\le n$, $k=0,1$, and
\begin{equation}\label{271eq}
\displaystyle{\sup_{r>0, Q\in\dd}\left\{ {1\over\oz(\drq)}\int_\trq G_0(y)\, {\a^2(y)\over\delta^2(y)}\, dy
\right\}^{1\over 2}\le\varepsilon_0,}
\end{equation}
then $\ol\in B_2(d\oz)$.
\end{thm}

We follow the ideas in~\cite{FeKP} (see also~\cite{K2} Theorem 2.7.1). Our proof is specialized to the case in which
the domain $D$ is the unit ball $B=B_1(0)$. Let $\W$ be as in (\ref{W}), $g\in C(\bb)$, and let $u_1$ be the solution
of
$$\left\{\begin{array}{rcl}
\Ll u_1&=&0\quad\text{ in }B\cr u_1 &=&g\quad\text{ on }\bb
\end{array}\right.$$
To prove Theorem~\ref{271} it is enough to show that for $\varepsilon_0$ small enough, there exist $C>0$ which depends
on $\alpha>0$ and the same parameters as $\varepsilon_0$ such that
\begin{equation}\label{cond274}
\| N_\alpha u_1\|_{L^2(\bb,d\oz)}\le C\,\| g\|_{L^2(\bb,d\oz)}.
\end{equation}
Indeed, (\ref{cond274}) is condition (ii) on Theorem~\ref{DP} when we take $\mu=\oz$ and $\omega=\ol$,
Theorem~\ref{271} then follows from Theorem~\ref{DP}.

Let now $u_0$ solve
$$\left\{\begin{array}{rcl}
\Lz u_0&=&0\quad\text{ in }B\cr u_0 &=&g\quad\text{ on }\bb
\end{array}\right.$$
Then from Corollary~\ref{grep} we have
$$u_1(x)=u_0(x)-\int_B G_0(x,y)\, \Lz u_1(y)\, dy = u_0(x)-F(x).$$

\begin{lem}\label{274} Under the hypothesis of Theorem~\ref{271} we have
$$N_\alpha^0 F(Q)\lesssim\varepsilon\, M_{\oz}(A_\alpha u_1)(Q),$$
where $N_\alpha^0$ denotes the modified nontangential maximal operator (Definition~\ref{enie}).
\end{lem}
\begin{lem}\label{275} Under the hypothesis of Theorem~\ref{271} we have
$$\int_{\bb}S_{c\alpha}^2 u_1\, d\oz\lesssim\, \int_{\bb}N_\alpha u_1^2\,d\oz.$$
\end{lem}

Let us take Lemmas~\ref{274} and~\ref{275} for granted, and let us show how we can then conclude the proof of
Theorem~\ref{271}. From Theorem~\ref{DP} we have
\begin{equation}\label{unot}
\| N_\alpha u_0\|_{L^p(\dd,d\oz)}\le C\, \| g\|_{L^p(\dd,d\oz)}.
\end{equation}
Now, from (\ref{unot}), Lemmas~\ref{enenie},~\ref{274} and~\ref{275}:
$$\begin{array}{rcl}
\displaystyle{\int_\bb (N_\alpha u_1)^2\,d\oz }&\lesssim &\displaystyle{\int_\bb (N^0_\alpha u_1)^2\,d\oz}\cr\cr
&\lesssim&\displaystyle{\int_\bb [(N^0_\alpha u_0)^2+(N_\alpha^0 F)^2]\,d\oz}\cr\cr &\lesssim&\displaystyle{\int_\bb
N_\alpha u_0^2\,d\oz+\,\varepsilon^2\,\int_\bb M_\oz(A_\alpha u_1)^2\,d\oz}\cr\cr &\lesssim&\displaystyle{\int_\bb
g^2\,d\oz+\,\varepsilon^2\,\int_\bb (A_\alpha u_1)^2\,d\oz}\cr\cr &\lesssim&\displaystyle{\int_\bb
g^2\,d\oz+\,\varepsilon^2\,\int_\bb (S_{c\alpha} u_1)^2\,d\oz}\cr\cr &\lesssim&\displaystyle{\int_\bb
g^2\,d\oz+\,\varepsilon^2\,\int_\bb (N_\alpha u_1)^2\,d\oz}
\end{array}$$
which proves Theorem~\ref{271} given that $\varepsilon$ is small enough. 

\subsection{Proof of Lemma~\ref{274}} 
 We assume, without lost of generality, that
\begin{equation}\label{assumed} \varepsilon(x)= 0\qquad\text{ for
}\delta(x)\ge\min\{ 1/2,r_0\}
\end{equation}
We note that from Lemma~\ref{lots}-(1) and Lemma~\ref{lem2} we have that if $\dist(y,\bb)=\delta(y)\le 1/4$ then
$${1\over \oz(\triangle_y)}{G_0(y)\over\delta(y)^2}\approx {\W(y)\over\W(B(y))}$$
Then, (\ref{271eq}) gives
\begin{equation}\label{note424}
\left\{ \int_{B(x)}\a^2(y){\W(y)\over\W(B(y))}\, dy\right\}^{1/2}\lesssim \varepsilon_0.
\end{equation}
Let $Q_0\in\bb$ and $x_0\in\Gamma_\alpha(Q_0)$. Denote by $\delta_0=\delta(x_0)$, $B_0=B_{\delta_0\over 6}(x_0)$,
$2B_0=B_{\delta_0\over 3}(x_0)$ and $B(x)=B_{\delta(x)\over 2}(x)$. Also, let $\tilde{G}(x,y)$ be the Green's function
for $\Lz$ on $B(x_0)=B_{\delta_0\over 2}$ and set
$$\begin{array}{lcl}
F_1(x) & = & \displaystyle{ \int_{2B_0} \tilde{G}(x,y)\,\Lz u_1(y)\, dy},\cr\cr F_2(x) & = & \displaystyle{ \int_{2B_0}
[G_0(x,y)-\tilde{G}(x,y)]\,\Lz u_1(y)\, dy},\cr\cr F_3(x) & = & \displaystyle{ \int_{B\backslash 2 B_0} G_0(x,y)\,\Lz
u_1(y)\, dy},
\end{array}$$
so that $F(x)=F_1(x)+F_2(x)+F_3(x)$, $x\in B(x_0)$. Remember $\varepsilon(x)=A_1(x)-A_0(x)$ and
$\a(x)=\sup_{B(x)}|\varepsilon(x)|$. If $x\in B_0$ we have
\begin{equation}\label{apsilon}
\begin{array}{rcl}
|\varepsilon(x)|&\le& \displaystyle{{C(n)\over |B(x)|}\, \int_{B(x)}\a(y)\, dy}\cr\cr &\le& \displaystyle{{C(n)\over
|B(x)|}\left\{\int_{B(x)}\a^2(y)\,{\W(y)\over \W(B(y))}\, dy\right\}^{1/2}\left\{\int_{B(x)}{\W(B(y))\over \W(y)}\,
dy\right\}^{1/2}}\cr\cr &\le&\displaystyle{{C\varepsilon_0}\left\{{\W(B(x))\over |B(x)|}\, {1\over |B(x)|}\int_{B(x)}
\W^{-1}(y)\, dy\right\}^{1/2}}\le\, C\, \varepsilon_0,
\end{array}\end{equation}
where we used (\ref{note424}), the doubling property of $\W$ and the fact that $\W\in A_2$. Note that $F_1$
verifies $\Lz F_1=\chi(2B_0)\, \Lz u_1$ in $B(x_0)$, $F_1=0$ on $\partial(B(x_0))$, with $\chi(2B_0)$ the
characteristic function of $2B_0$. From the weighted Sobolev inequality (Theorem 1.2 in~\cite{FKS}),
Theorem~\ref{strrong} and (\ref{apsilon}) we have
\begin{equation}\label{F1}
\begin{array}{rcl}
\left\{\int_{2B_0} F_1^2(x)\, {\ww}\, dx\right\}^{1\over 2} &\lesssim& \delta_0\,\left\{\int_{B(x_0)}
|\nabla F_1(x)|^2\, {\ww}\, dx\right\}^{1\over 2} \cr\cr &\lesssim& \left\{\int_{B(x_0)} |\Lz
F_1(x)|^2\,\delta^4(x)\, {\ww}\, dx\right\}^{1\over 2} \cr\cr &\lesssim& \left\{\int_{2B_0}
\varepsilon(x)^2 |\nabla^2 u_1(x)|^2\,\delta^4(x)\, {\ww}\, dx\right\}^{1\over 2}\cr\cr &\lesssim&
\varepsilon_0\, A_\alpha u_1(Q_0).
\end{array}\end{equation}

Let $v_x(y)=G_0(x,y)-\tilde{G}(x,y)$, $x,\, y\in 2B(x_0)$. Then if $\Lz^*$ denotes the adjoint operator to $\Lz$, we
have $\Lz^* v_x = 0$ in $B(x_0)$, and $v_x\ge 0$ in $B(x_0)$. In particular, $\tilde{v}(y)=v_x(y)/\W(y)$ is a
nonnegative {\em nor\-ma\-li\-zed adjoint solution} (n.a.s.) of $\Lz$ in $B(x_0)$ (see Definition~\ref{nas}). From the
maximum principle for n.a.s., Harnack inequality for n.a.s. and Lemma~\ref{lem2} we have
\begin{equation}\label{mhnas}
\tilde{v}\le \max_{y\in\, 2B_0}\tilde{v}\le C\, {G_0(x,\overline{y})\over\W(\overline{y})}\approx
{\omega^x(\triangle_{\overline{y}})\,\delta^2(\overline{y})\over\W(B(\overline{y}))},
\end{equation}
where $\overline{y}\in\partial (2B_0)$ verifies $\delta(\overline{y})={\rm dist}(\partial(2B_0),\bb)$. From
Lemma~\ref{lots}-(2) and Harnack inequality it follows
$\omega^x(\triangle_{\overline{y}})\approx 1$. Then, from Lemma~\ref{lem2}, (\ref{mhnas}) and (\ref{apsilon}) we get
for $x\in B_0$
\begin{equation}\label{F2}\begin{array}{rcl}
|F_2(x)|&\lesssim& \displaystyle{\varepsilon_0\int_{2B_0}\delta(\overline{y})^2\,
\omega^x(\triangle_{\overline{y}})\,|\nabla^2u_1(y)|\,{\W(y)\over\W(B(\overline{y}))}\, dy}\cr\cr &\lesssim&
\displaystyle{\varepsilon_0\left\{\int_{2B_0}\delta({y})^4\,|\nabla^2u_1(y)|^2\, {\W(y)\over\W(B(y))}\, dy\right\}^{1\over
2}\left\{\int_{B_0}{\W(y)\over\W(B(y))}\, dy\right\}^{1\over 2}}\cr\cr &\lesssim&\displaystyle{\varepsilon_0\, A_\alpha
u_1(Q_0)},
\end{array}\end{equation}
Now define
$$\O_0=\left( B\bigcap B_{\delta_0\over 2}(Q_0)\right),\quad \O_j=\left( B\bigcap B_{2^{j-1}\delta_0}(Q_0)\right)\backslash\left( 2B_0\bigcup B_{2^{j-2}\delta_0}(Q_0)\right),$$
$j=1,\, 2,\,\cdots,\, j_0-1$, and
$$\O_{j_0}=B\backslash B_{2^{j_0-2}\delta_0}(Q_0),$$
where ${1\over 2}<\delta_0\, 2^{j_0-1}\le 1$.  We set  
$$F_3^j(x)=\int_{\Omega_j}G_0(x,y)\, \Lz u_1(y)\, dy,\quad j=0,\,1,\,\cdots,\, j_0.$$
Thus
\begin{equation}\label{Fsum}
F_3(x)=\sum_{j=0}^{j_0} F_3^j(x).
\end{equation}
We now need to define the notion of a dyadic grid:
\begin{defn}\label{dyadic} A {\em{ dyadic grid}} on $\b$ is a collection of Borel sets $\I=\bigcup_{k=1}^\infty\I_k$, such that
\begin{itemize}
  \item[{\rm{(i)}}] $\b=\bigcup_{I\in\I_k}I$,\quad $k=1,\, 2,\, \cdots$,
  \item[{\rm{(ii)}}] $I,\, J\in\I_k$, $I\neq J$ $\Rightarrow$ $I$ and $J$ have disjoint interiors, 
  \item[{\rm{(iii)}}] $I\in\I_k\Rightarrow \diam(I)\approx 2^{-k}$,
  \item[{\rm{(iv)}}] for every $I\in \I_k$, $\exists !$ $I'\in \I_{k-1}:$ $I\subset I'$, $k=2,\,3,\,\cdots,$
  \item[{\rm{(v)}}] for every $I\in\I$ there exists a boundary ball $\triangle_I$ such that
$$\triangle_I\subset I\subset c\,\triangle_I.$$
where $c>1$ is a universal constant.  
\end{itemize}
If $I\in\I$ we say that $I$ is dyadic.
\end{defn}
We consider now a Whitney decomposition of $B$ into cubes, $B=\bigcup_{Q\in\Q} Q$, as in the proof of Theorem~\ref{whit}, it is easy to see that there exists a dyadic grid on $\b$ such that for every $I\in\I$ we can assign a cube $I^+\in\Q$ so that
\begin{equation}\label{ring}
\bigcup_{k\ge 1}\bigcup_{I\in \I_k} I^+ = B\backslash B_{1\over 2}(0).
\end{equation}
Moreover, the correspondence $I\mapsto I^+$ can be defined so that for any boundary ball $\drq$ we have $\trq\subset\bigcup_{I\in\I,\,I\subset\drq}I^+$. 
Remember that by (\ref{assumed}) $\varepsilon\equiv 0$ in $B_{1\over 2}(0)$. Let $\tz=\triangle_{\delta_0\over
2}(Q_0)=B_{\delta_0\over 2}\bigcap\bb$, and suppose $I\subset\tz$ dyadic. For $x\in B_0$ and  $y\in I^+$ we have from
Lemma~\ref{lots}-(1) and Lemma~\ref{lem2}
\begin{equation}\label{equivalences}
\displaystyle{G_0(x,y)\approx \delta^2(y)\,
{\oz(\triangle_y)\over \oz(\tz)}{\W(y)\over\W(B(y))}\approx {G_0(y)\over \oz(\tz)}}.\end{equation} To estimate
$F_3^0(x)$ we will use a ``stopping time'' argument. For $j=0,\,\pm 1,\,\pm 2,\,\cdots$, let
$$\begin{array}{ll}
\oo_j=\{ Q\in\tz\,:\, A_\alpha u_1(Q)> 2^j\},&\cr\cr \tilde{\oo}_j=\{ Q\in\tz\,:\, M_\oz(\chi(\oo_j))(Q)> c^*\},&\cr\cr
\J_j=\{ I \text{ dyadic},\,I\subset\tz \,:\, \oz(I\bigcap\oo_j)\ge{1\over 2}\,\oz(I)\text{ but
}\oz(I\bigcap\oo_{j+1})<{1\over 2}\,\oz(I)\}.&
\end{array}$$
If $I\in\J_j$, by Definition~(\ref{dyadic})-(v) and the doubling property of $\oz$ it easily follows that for any $Q\in I$
$$M_\oz(\chi(\oo_j))(Q)>c^*$$
for some $c^*>0$ depending only on dimension and ellipticity. Therefore, ta\-king this choice of $c^*$ in the definition of
$\tilde{\oo_j}$, we have $I\subset\tilde{\oo}_j,\text{ whenever }I\in\J_j.$ 
Now, from  (\ref{equivalences}) we get
$$\begin{array}{rl}
|F_3^0(x)|& =\displaystyle{\left|\sum_{\stackrel{\displaystyle{I\subset\tz}}{I\text{
dyadic}}}\int_{I^+\bigcap\Omega_0}G_0(x,y)\, \Lz u_1(y)\, dy\right|}\cr\cr
\lesssim &
{1\over\oz(\tz)}\,\displaystyle{{\sum_j\sum_{I\in\J_j}}\oz(I)\int_{I^+\bigcap\Omega_0}
\delta^2(y)\,\varepsilon(y)\,|\nabla^2u_1(y)|\,{\W(y)\over\W(B(y))}\,dy}.
\end{array}$$

On the other hand, if $I\subset\tz$ is dyadic then $I\in\J_j$ for some $j$. Now we set
$$\overline{\J_j}=\{J\text{ dyadic: }J\in\J_j,\text{ and }J\subset  J',\, J'\text{ dyadic, }J'\in\J_j,\Rightarrow J'=J\},$$
that is, $\overline{\J_j}$ is the collection of ``maximal'' dyadic sets in $\J_j$. It is clear that $\J_j=\bigcup_{J\in\overline{\J_j}} J$, with a disjoint union. Let $T(J)=\bigcup_{\stackrel{\displaystyle{I\subset
J}}{I\text{ dyadic}}}I^+$, where $J$ is dyadic, then if $J\in\overline{\J_j}$, from (\ref{271eq})  and (\ref{equivalences}) we have
\begin{equation}\label{leps}
\displaystyle{\left\{\sum_{\stackrel{\displaystyle{I\text{ dyadic}}} {I\subset
J}}\oz(I)\,\int_{I^+\bigcap\Omega_0}\a^2(y)\,{\W(y)\over\W(B(y))}\,dy\right\}^{1\over
2}\lesssim\left\{\int_{T(J)\bigcap\Omega_0}G_0\,{\a^2\over\delta^2}\,dy\right\}^{1\over 2}\lesssim
\varepsilon_0\, \left(\oz(J)\right)^{1/2}.}\end{equation}
By H\"older inequality and (\ref{leps}),
$$\begin{array}{rcl}
|F_3^0(x)|& \lesssim &
\displaystyle{{\varepsilon_0\over\oz(\tz)}\,\sum_j\sum_{J\in\overline{\J_j}}\left(\oz(J)\right)^{1/2}\left\{\sum_{J\supset I\in \J_j} \oz(I)\int_{I^+\bigcap\Omega_0} \delta^4(y)\,|\nabla^2u_1(y)|^2\,{\W(y)\over\W(B(y))}\,dy\right\}^{1/2}}.
\end{array}$$
Now we make the following observation:
There exists $\alpha=\alpha(n)$ such that if $I\subset\bb$ is dyadic, and $E\subset I$ with $2\,\oz(E)\ge\oz(I)$, then
$$\int_{I^+} f(y)\,\oz(\triangle_y)\, dy\lesssim\int_E\int_\Gaq f(y)\, dy\, d\oz(Q).$$
This is a consequence of Fubini's theorem and the fact that for appropriate $\alpha=\alpha(n)$ we have $I^+\subset\Gaq$
for all $Q\in I$.
 From the weak type inequality of the maximal operator $M_\oz$ we have $\oz(\tilde{\oo_j}\backslash\oo_{j+1})\lesssim\oz(\oo_j)$. Then, since $\oz(I\backslash\oo_{j+1})\ge{1\over 2}\,\oz(I)$ and $I\subset\tilde{\oo}_j$ for all $I\in\J_j$ we have
\begin{equation}\label{F30}\begin{array}{rcl}
|F_3^0(x)|&\lesssim &\displaystyle{ \sum_j{\varepsilon_0\,\oz(\oo_j)^{1\over
2}\over\oz(\tz)}\,\left\{\int_{\tilde{\oo_j}\backslash\oo_{j+1}}\int_\Gaq|\nabla^2u_1(y)|^2\, {\delta^4(y)\,\W(y)\over\W(B(y))} \,dy\, d\oz\right\}^{1\over 2}}\cr\cr 
&\lesssim &\displaystyle{ {\varepsilon_0\over\oz(\tz)}\,\sum_j\oz(\oo_j)^{1\over 2}\,\left\{\int_{\tilde{\oo_j}\backslash\oo_{j+1}}A^2_\alpha u_1(Q)\, d\oz(Q) \right\}^{1\over 2}}\cr\cr 
&\lesssim&\displaystyle{ {\varepsilon_0\over\oz(\tz)}\,\sum_j\oz(\oo_j)\, 2^j}\cr\cr &\lesssim &\displaystyle{
{\varepsilon_0\over\oz(\tz)}\int_\tz A_\alpha u_1(Q)\, d\oz(Q)}\lesssim \varepsilon_0\, M_\oz( A_\alpha u_1 )(Q_0).
\end{array}\end{equation}
Now we claim that for some $\theta>0$ we have
\begin{equation}\label{F3j}
|F_3^j(x)|\lesssim 2^{-j\theta}\,\varepsilon_0\, M_\oz(A_\alpha u_1)(Q_0),\qquad j=1,\, 2\, \cdots .
\end{equation}
If we take (\ref{F3j}) for granted, adding in $j$ in (\ref{F3j}) and from (\ref{F30}) and (\ref{Fsum}) we obtain
\begin{equation}\label{F3}
|F_3(x)|\lesssim\varepsilon_0\, M_\oz(A_\alpha u_1)(Q_0),\qquad x\in B_0.\end{equation}
Then, from $F(x)=F_1(x)+F_2(x)+F_3(x)$, $x\in B_0$, (\ref{F1}), (\ref{F2}), (\ref{F3}) and the doubling property of $\W$ it follows that
$$N_\alpha^0 F(Q_0)\lesssim\varepsilon_0\, M_\oz(A_\alpha u_1)(Q_0),\qquad x\in B_0,$$
which proves Lemma~\ref{274}.

To show (\ref{F3j}), we proceed as in the estimate of $F_3^0$. We set
$\tj=\triangle_{2^{j-1}\delta_0}(Q_0)=B_{2^{j-1}\delta_0}(Q_0)\bigcap\bb$, and $\tjz=\tj\backslash\triangle_{j-1}$.
From Definition~\ref{dyadic} and simple geometrical considerations it follows that there exists $\alpha>0$ such that
$$\O_j\subset(\Gamma_\alpha(Q_0)\bigcap\O_j)\bigcup(\bigcup_{\stackrel{
\displaystyle{I\subset\tjz}}{I\text{ dyadic}}}I^+).$$ Let $x_j\in (\Gamma_\alpha(Q_0)\bigcap\O_j)=\O_j^0$, $j=1,\, 2,\,
\cdots$, since $G_0(\cdot,y)$ is a nonnegative solution of $\Lz u=0$ in $B\backslash\{ y\}$ vanishing on $\b$, $G(\cdot,y)$ is H\"older continuous up to the boundary~\cite{Ba}, moreover, we have
$$G_0(x,y)\lesssim 2^{-j\theta}\, G_0(x_{j-1},y),\qquad y\in \O_j^0,$$
where $C$ and $\theta$ only depend on ellipticity and dimension. From this inequality and Harnack inequality for
nonnegative solutions we have $G_0(x,y)\lesssim 2^{-j\theta}\, G_0(x_{j+1},y)$. From Lemma~\ref{lem2} we get
$$G_0(x_{j+1},y)\approx \oz^{x_{j+1}}(\tj)\,\W(y)\,\W(B(y))^{-1}\delta(y)^2 \approx {\oz^{x_{j+1}}(\tj)\over \oz(\tj)}\,G_0(y).$$
From Lemma~\ref{lots}-(2) we have $\oz^{x_{j+1}}(\tj)\approx 1$, thus
\begin{equation}\label{thetaG}
G_0(x,y)\lesssim{2^{-j\theta}\over \oz(\tj)}\, G_0(y),\qquad x\in B_0,\  y\in\O_j^0.
\end{equation}
We have for $x\in B_0$
\begin{equation}\label{F3j0}
\begin{array}{rcl}
\displaystyle{|\int_{\O_j^0}\Lz u_1(x)\, G_0(x,y)\, dy|}&\lesssim&\displaystyle{
{2^{-j\theta}\over\oz(\tj)}\int_{\O_j^0}\varepsilon(x)\,|\nabla^2u_1(x)|\, G_0(y)\, dy}\cr\cr &\lesssim&
\displaystyle{{2^{-j\theta}\over\oz(\tj)}\left\{ \int_{\O_j^0}|\nabla^2u_1(x)|\, \delta(y)^4\,{\W(y)\over\W(B(y))}\,
dy\right\}^{1\over 2}}\cr &&\qquad\displaystyle{\left\{\int_{\O_j^0} \a^2(x)\,\delta(y)^{-4}\, {G_0}^2(y)\,
{\W(B(y))\over\W(y)}\,dy\right\}^{1\over 2}}\cr\cr &\lesssim& \displaystyle{{2^{-j\theta}}\left\{
\int_{\O_j^0}|\nabla^2u_1(x)|\, \delta(y)^4\,{\W(y)\over\W(B(y))}\, dy\right\}^{1\over 2}}\cr &&\qquad\displaystyle{\left\{
{1\over\oz(\tj)}\int_{\O_j^0} \a^2(x)\,\delta(y)^{-2}\, {G_0}(y) \,dy\right\}^{1\over 2}}\cr\cr &\lesssim&
\displaystyle{\varepsilon_0 \,{2^{-j\theta}}\, A_\alpha(u_1)(Q_0).}
\end{array}
\end{equation}
On the other hand, an argument similar to the one applied to obtain the bound for $F_0$ and  a consideration in the
spirit of (\ref{thetaG}) yields
\begin{equation}\label{F3j1}
\displaystyle{|\int_{\O_j\backslash\Gamma_\alpha(Q_0)}\Lz u_1(x)\, G_0(x,y)\, dy|\lesssim \varepsilon_0\,
{2^{-j\theta}}\, M_\oz(A_\alpha(u_1))(Q_0).}
\end{equation}
(\ref{F3j}) follows from (\ref{F3j0}) and (\ref{F3j1}), this concludes the proof of Lemma~\ref{274}.

\subsection{Proof of Lemma~\ref{275}} From the identity $\Lz(u_1^2)=2\,A_0\,\nabla u_1\cdot\nabla u_1+2\, u_1\,\Lz u_1,$ we have
$$|\nabla u_1|^2\lesssim 2\,A_0\,\nabla u_1\cdot\nabla u_1=\Lz(u_1^2)-2\, u_1\,\Lz u_1,$$
We let $B^*=B_{1-r_0}$, by Fubini's theorem and the fact that $\varepsilon(x)=0$ in $B_{1\over 2}$ we get
$$\begin{array}{rcl}
\displaystyle{\int_\bb S_{c\alpha}^2u_1\, d\oz}&\lesssim&\displaystyle{\int_{B\backslash B^*} |\nabla
u_1|^2\,\delta(x)^2\,{\ww}\,\oz(\triangle_x)\, dx}\cr\cr &\lesssim&\displaystyle{\int_{B\backslash B^*} \{ \Lz(u_1^2)-2\,
u_1\,\Lz u_1\}\,\delta(x)^2\,{\ww}\,\oz(\triangle_x)\, dx}\cr\cr &\lesssim&\varepsilon_0\,\displaystyle{\int_{B\backslash
B_{1\over 2}} |u_1|\,|\nabla^2 u_1| \,\delta(x)^2\,{\ww}\,\oz(\triangle_x)\, dx}
\end{array}$$
where we used 
Lemma~\ref{lem2} and $\int_B
\Lz(u_1^2)\,G_0(x)\, dx\le 0$. Now we apply a ``stopping time'' argument as in the proof of (\ref{F30}) to obtain
$$\label{stop2}\displaystyle{\int_\bb S_{c\alpha}^2u_1\, d\oz}\lesssim \varepsilon_0\,\displaystyle{\int_\bb N_\alpha
u_1\cdot A_\alpha u_1\,d\oz}.$$
From Theorem~\ref{whit} and the inequality $|a\, b|\le \mu^{-1}\,a^2+\mu\, b^2$, $\mu>0$
we get
$$\begin{array}{rcl}
\displaystyle{\int_\bb S_{c\alpha}^2u_1\, d\oz}&\le&C\,\varepsilon_0\,\displaystyle{\int_\bb N_\alpha u_1\cdot S_{c\alpha}
u_1\,d\oz}\cr\cr &\le&\displaystyle{C\,{\varepsilon_0\over\mu}\,\int_\bb (N_\alpha u_1)^2
\,d\oz+C\,\varepsilon_0\,\mu\,\int_\bb S^2_{c\alpha} u_1\,d\oz},
\end{array}$$
the Lemma follows from choosing $\mu$ so that $C\,\varepsilon_0\,\mu={1\over 2}$.

\section{Proof of Theorem~\ref{main}}
We will obtain Theorem~\ref{main} as a consequence of the following special case in the spirit of~\cite{RF}.
\begin{thm}\label{mainrf} Suppose that $\oz\in A_\infty(d\sigma)$ and let $E(Q)$ be given by
$$E_r(Q)=\left\{ \int_{\Gaq\bigcap B_r(Q)} {\a^2(x)\over \delta^n(x)}\,
dx\right\}^{1\over 2}\qquad r>0,\, Q\in\dd.$$ There exists $\varrho_0=\varrho_0(n,\lambda,\Lambda,D,\zeta,\kappa)>0$
such that if $a_k^{i,j}\in\bmss(\r)$, $1\le i,j\le n$, $k=0,1$, and
\begin{equation}\label{carlcone}
\displaystyle{\sup_{Q\in\dd}E_{r_0}(Q)=M_1<\infty,}
\end{equation}
then $\ol\in A_\infty(d\sigma)$. Here, $\zeta$ and $\kappa$ are the $A_\infty$ constants of $\oz$ with respect to
$\sigma$ as given in (\ref{ainf}).
\end{thm}
We postpone the proof of Theorem~\ref{mainrf} to the next section, and we now show how Theorem~\ref{main} follows from
this result. We assume that $\a(x)\equiv 0$ if $\delta(x)>r_0/2$. By Lemma~\ref{wlog}, this assumption does not bring
any loss of generality.  Fix $Q\in\dd$ and let $0<r\le r_0/2$, by Fubini's theorem and (\ref{carl}) we have
$$\begin{array}{rcl}
\displaystyle{{1\over\sigma(\drq)}\int_\drq E_r^2(P)\, d\sigma(P)}&\lesssim&
\displaystyle{{1\over\sigma(\drq)}\int_\ttr {\a^2(x)\over\delta(x)}\, dx}\lesssim M^2.
\end{array}$$
So there exists a closed set $F\subset\drq$ such that $2\,\sigma(F)>\sigma(\drq)$ and $E(P)\lesssim M$ for all $P\in
F$. Now, we need to introduce a ``saw-tooth'' region $\O=\O(F,r)$ over $F$, that is, for given $0<\alpha<\beta$, $\O$ verifies (see~\cite{DJK},~\cite{EK}):
\begin{enumerate}%
  \item[(i)] for suitable $\alpha'$, $\alpha''$, $c_1$, $c_2$ with $\alpha<\alpha'<\alpha''<\beta$
$$\bigcup_{P\in E} \{ \overline{\Gamma}_{\alpha'}(P)\bigcap B_{c_1\, r}(P)\}\subset\O\subset \bigcup_{P\in E}\{
\overline{\Gamma}_{\alpha''}(P)\bigcap B_{c_2\, r}(P)\};$$
  \item[(ii)] $\partial\O\bigcap\dd=F$;
  \item[(iii)] there exists $x_0\in\O$ with $\dist(x_0,\partial\O)\approx r$;
we call $x_0$ the center of $\O$; 
  \item[(vi)] $\O$ is a Lipschitz domain with Lipschitz constant which depends only on $B$.
\end{enumerate}
\noindent Let $\tilde{\L}=\sum_{i,j=1}^n {\tilde{a}^{i,j}}\pij $ where $\tilde{A}=\{ \tilde{a}^{i,j}(x)\}_{i,j=1}^n$ and
$$\tilde{A}(x)=\left\{\begin{array}{ll}
\displaystyle{A_1(x)\quad }&\displaystyle{x\in\O(F,r)}\cr
\displaystyle{A_0(x)\quad }&\displaystyle{x\in B\backslash\O(F,r)}
\end{array}\right.$$
From the definition of $\a(x)$ it is easy to see that
\begin{equation}\label{proportion}
\displaystyle{|\{ y\in B(x):\a(y)\ge\a(x)\}|\gtrsim|B(x)|.}
\end{equation}
From (\ref{proportion}) and (\ref{carl}) we have that for all $x\in D$
$$\a^2(x)\lesssim {1\over |B(x)|}\int_{B(x)}\a^2(y)\, dy\lesssim M^2.$$
We set $\tilde{\a}(x)=\sup_{y\in B(x)}|\tilde{A}(y)-A_0(y)|$, then $\tilde{\a}(x)\le \a(x)\lesssim M$ for all $x\in D$.
Note that since for all $1\le i,j\le n$ we have
$$\begin{array}{rcl}
\tilde{a}^{i,j}&=&a_0^{i,j}+(a_1^{i,j}-a_0^{i,j})\chi_{\O(F,r)}
\end{array}$$
where $\chi_{\O(F,r)}$ is the characteristic function of $\O(F,r)$, we have $\eta_{\tilde{a}}\le\eta_{a_0}+C\, M$,
where $\eta_{\tilde{a}}$ and $\eta_{a_0}$ are the $\bmo$ moduli of continuity of $\tilde{A}$ and $A_0$, respectively.
Therefore, if $M$ is small enough, say $M\le C^{-1}\varrho^*$, we have that the coefficients of $\tilde{\L}$ belong to
the space $\bmss$.

Our saw-tooth region $\O(F,r)$ can be constructed so that for any $Q\in\dd$ such that
$x\in\Gaq\bigcap\O(F,r)\neq\emptyset$, there exists $P\in F$ such that $B(x)\subset\Gamma_\alpha(P)$. From this
observation, we have that if $\tilde{E_r}(P)$ is as in the definition of $E_r(P)$ above but replacing $a$ by
$\tilde{a}$, then $\tilde{E_r}(P)\lesssim M$ for all $P\in\dd$. Therefore, from Theorem~\ref{mainrf} we have that the
harmonic measure $\tilde{\omega}$ for $\tilde{\L}$ on $\dd$ is in $A_\infty(d\sigma)$ and by a known result in the
theory of weights~\cite{Mu} there exists $\theta>0$ and $c>0$ (depending only on $\zeta$ and $\kappa$) such that
\begin{equation}\label{ainfp}
\left({\tilde{\omega}(Z)\over\tilde{\omega}(\drq)}\right)^\theta\ge c\, {\sigma(Z)\over \sigma(\drq)}
\end{equation}
for any set $Z\subset\drq$. Since ${\sigma(F)\over\sigma(\drq)}>{1\over 2}$ we have that if
${\sigma(E)\over\sigma(\drq)}>{3\over 4}$ then ${\sigma(E\bigcap F )\over\sigma(\drq)}>{1\over 4}$. Therefore, from
(\ref{ainfp}) we get
$${\tilde{\omega}(E\bigcap F)\over\tilde{\omega}(\drq)}\gtrsim 1.$$
Let $x_0$ be the ``center'' of the saw-tooth region $\O(R,r)$ and denote by
$\tilde{\omega}_\O$ the harmonic measure for $\tilde{\L}$ on $\partial\O(F,r)$ evaluated at $x_0$. From the ``main lemma'' in ~\cite{EK} (see also~\cite{DJK}) we have
$${\tilde{\omega}(E\bigcap F)\over \tilde{\omega}(\drq)}\le C\, \left(\tilde{\omega}_\O(E\bigcap F)\right)^{\vartheta}.$$
Thus, $\tilde{\omega}_\O(E\bigcap F)\ge C$. Since $\Ll=\tilde{\L}$ in $\O(F,r)$, we obtain $\ol_\O(E\bigcap F)\ge C$,
where $\ol_\O$ denotes the harmonic measure on $\partial\O(F,r)$ for the operator $\Ll$. From the scaling of the
harmonic measure (Lemma~\ref{lots}-(1)) and the maximum principle we get
$${\ol(E\bigcap F)\over \ol(\drq)}\gtrsim \ol^{x_0}(E\bigcap F)\ge \ol_\O(E\bigcap F)\gtrsim 1$$
and then:
$${\ol(E)\over\ol(\drq)}\ge \kappa_0,$$
for some positive $\kappa_0$. This shows that condition (\ref{ainf}) holds for the measures $\ol$ and $\sigma$ with
$\zeta={3\over 4}$ and $\kappa=\kappa_0$. Therefore, $\ol\in A_\infty(d\sigma)$ as wanted, in the case $M\le
C^{-1}\varrho^*$.

For the general case, define $\L_t=(1-t)\, \Lz+t\, \Ll$ for $0\le t\le 1$, let $K$ be a positive integer such that
${K^{-1}}\, M\le C^{-1}\varrho^*$ and for integers $0\le l < K$ let $\a_l(x)=\sup_{B(x)}|A_{{l+1}\over K}(y)-A_{l\over
K}(y)|$, where $A_t=(1-t)\, A_0+t\, A_1$. Note that since the set $\bms$ is convex we have that $A_t\in\bms$, $0\le
t\le 1$, by the same consideration, the matrices $A_t$ are uniformly elliptic with ellipticity $\lambda$ and their
entries are bounded by $\Lambda$. Denote by $\omega_t$ the harmonic measure in $D$ with respect to $\L_t$, $0\le t\le
1$. Then since
$$\a_l(x)=\sup_{y\in B(x)}|{1\over K}\,(A_1(y)-A_0(y))|\le {1\over K}\, \a(x),\qquad 0\le l<K,$$
we have that the pairs of operators $\L_{l\over K}$ and $\L_{l+1\over K}$, $0\le l<K$ verify
$$\sup_{\stackrel{0<r<r_0}{Q\in\dd}}h^l(r,Q)=M/K\le C^{-1}\varrho^*<+\infty,$$ where
$$h^l(r,Q)=\left\{ {1\over{\sigma(\drq)}}\int_{\trq} {\a_l^2(x)\over\delta(x)}\, dx\right\}^{1/2},\qquad 0\le l<K,$$
from the previous special case we have
$$\omega_{l\over K}\in A_\infty(d\sigma)\Rightarrow\omega_{{l+1}\over K}\in
A_\infty(d\sigma),\qquad l=0,\, ,1,\, \cdots\, K-1.$$ So $\ol\in A_\infty(d\sigma)$, this finishes the proof of
Theorem~\ref{main}. 

\section{Proof of Theorem~\ref{mainrf}}

Now we will show that Theorem~\ref{mainrf} follows from Theorem~\ref{271}. It is clear from the proof of
Theorem~\ref{271} (see Section~\ref{proof271}) that we might replace $\a(x)$ in (\ref{271eq}) by $\a_0(x)=\sup_{y\in
B_0(x)}|A_0(y)-A_1(y)|$ where $B_0(x)$ is the ball centered at $x$ with radius $\delta(x)/c$ for any fixed constant
$c\ge 2$. We claim that if we take $c=8$ in the definition of $\a_0$, $M_1$ is given in by (\ref{carlcone}) and $0<r\le
r_0$ we have
\begin{equation}\label{weightineq}
\displaystyle{ {1\over \oz(\drq)}\int_{\trq}\a_0^2(x){G_0(x)\over\delta(x)^2}\, dx\lesssim M_1^2.}
\end{equation}
In fact, from Lemma~\ref{lem2} we have
$$\begin{array}{rcl}
\displaystyle{\int_{\trq}\a_0^2(x){G_0(x)\over\delta(x)^2}\, dx}&\lesssim &
\displaystyle{\int_{\trq}\a_0^2(x)\,{\ww}\,\oz(\triangle_x) dx.}
\end{array}$$
From Besicovitch's covering lemma (I.8.17 in~\cite{St2}), we can find a sequence $\{x_j \}_{j=1}^\infty\subset\trq$
such that $\trq\subset\bigcup_j B_j$, where $B_j=B_{\delta(x_j)\over 8}(x_j)$ and the balls $B_j$ have finite
overlapping. Then
$$\begin{array}{rcl}
\displaystyle{\int_{\trq}\a_0^2(x){G_0(x)\over\delta(x)^2}\, dx}&\lesssim &
\displaystyle{\sum_j\int_{B_j}\a_0^2(x)\,{\ww}\,\oz(\triangle_x) dx}\cr\cr 
&\lesssim&\displaystyle{\sum_j\,{\oz(\triangle_{x_j})\over\W(B(x_j))}\int_{B_j}\a_0^2(x)\,{\W(x)} dx}\cr\cr 
&\lesssim &\displaystyle{\sum_j\,{\oz(\triangle_{x_j})}\,\a^2({x}_j)}\cr\cr 
&\lesssim &\displaystyle{\sum_j\,{\oz(\triangle_{x_j})}\int_{B({x}_j)}{\a^2(x)\over \delta(x)^n }\, dx}\cr\cr 
&\lesssim &\displaystyle{\int_{\drq}E_r^2(P)\, d\oz(P),}
\end{array}$$
where we have used Fubini's theorem; (\ref{weightineq}) now follows from (\ref{carlcone}). Then, if $M_1$ is small enough, say $M_1\le C^{-1}\varepsilon_0$,
from Theorem~\ref{271} we have that $\ol\in B_2(d\oz)$ and since $\oz\in A_\infty(d\sigma)$, we have $\ol\in
A_\infty(d\sigma)$, proving Theorem~\ref{mainrf} in the case $M_1\le C^{-1}\varepsilon_0$.

For the general case, we define as before $\L_t=(1-t)\, \Lz+t\, \Ll$ for $0\le t\le 1$, let $K$ be a positive integer
such that ${C\, K^{-1}}\, M_1\le\varepsilon_0$ and for integers $0\le l < K$ let $\a_l(x)=\sup_{B(x)}|A_{{l+1}\over
K}(y)-A_{l\over K}(y)|$, where $A_t=(1-t)\, A_0+t\, A_1$. Denote by $\omega_t$ the harmonic measure in $D$ with respect
to $\L_t$, $0\le t\le 1$. Then since
$$\a_l(x)=\sup_{y\in B(x)}|{1\over K}\,(A_1(y)-A_0(y))|\le {1\over K}\, \a(x),$$
from the previous result we have
$$\omega_{l\over K}\in A_\infty(d\sigma)\Rightarrow\omega_{{l+1}\over K}\in
A_\infty(d\sigma),\qquad l=0,\, 1,\, \cdots,\, K-1.$$ So $\ol\in A_\infty(d\sigma)$, this completes the proof of
Theorem~\ref{mainrf}. 

\subsection*{Acknowledgment}
The results of this paper were obtained during my Ph.D. stu\-dies at
the University of Minnesota and are also contained in my
thesis~\cite{Cr}. I would like to express
immense gratitude to my supervisor Carlos Kenig, from the University of Chicago, whose
guidance and support were essential for the successful completion of
this project.

\end{document}